\documentclass[12pt]{amsart}

\newcommand{\Y}{{\mathcal Y}}
\newcommand{\LL}{{\mathcal L}}
\newcommand{\MM}{{\mathcal M}}

\newcommand{\Abar}{{\overline{A}}}

\newcommand{\pibar}{{\overline{\pi}}}
\newcommand{\Diff}{\operatorname{Diff}}
\newcommand{\D}{{\mathcal D}}
\newcommand{\Dbar}{{\overline{\D}}}
\newcommand{\fancyS}{{\mathcal S}}

\newcommand{\C}{{\mathbf C}}

\newcommand{\kbar}{{\overline{k}}}

\newcommand{\Z}{{\mathbf Z}}
\newcommand{\R}{{\mathbf R}}
\newcommand{\PP}{{\mathbf P}}

\newcommand{\G}{{\mathbf G}}
\newcommand{\tors}{{\operatorname{tors}}}

\newcommand{\Gal}{\operatorname{Gal}}

\newcommand{\nichts}{{\left.\right.}}
\newcommand{\isom}{\cong}
\newcommand{\tensor}{\otimes}

\newtheorem{theorem}{Theorem}
\newtheorem{lemma}[theorem]{Lemma}
\newtheorem{cor}[theorem]{Corollary}

\newtheorem{conj}[theorem]{Conjecture}

\theoremstyle{definition}

\theoremstyle{remark}
\newtheorem{rem}{Remark$\!\!$}	
\newtheorem{rems}{Remarks$\!\!$}	

\begin{document}

\title{Mordell-Lang plus Bogomolov II: the division group}
\author{Bjorn Poonen}
\thanks{This research was supported
by NSF grant DMS-9801104
and by a Sloan Research Fellowship.}
\address{Department of Mathematics, University of California,
		Berkeley, CA 94720-3840, USA}
\email{poonen@math.berkeley.edu}
\date{August 23, 1998}

\maketitle

\section{Introduction}

In~\cite{mlb}, we formulated a conjecture for semiabelian varieties
which would include both the Mordell-Lang conjecture and the
Bogomolov conjecture,
and we proved this conjecture in the case where
the semiabelian variety was {\em almost split}
(i.e., isogenous to the product of an abelian variety and a torus),
and where we had only a finitely generated subgroup
instead of its division group.

In this paper we prove the conjecture for the full division group.
The ``almost split'' hypothesis remains,
but this is only because
it is under this hypothesis that the Bogomolov conjecture
and an equidistribution theorem
have been proved so far~\cite{chambertloir2}.
In fact, if we assume the Bogomolov conjecture
and an equidistribution theorem for general semiabelian varieties,
then we can prove our conjecture entirely.

Our proof requires also the Mordell-Lang conjecture
(which is completely proven).
But we do not need its full strength:
we use only the ``Mordellic'' (finitely generated) part.
Hence our proof gives a new reduction of the full Mordell-Lang conjecture
to the Mordellic part,
at least for almost split semiabelian varieties.

We postpone the precise statement of our result
until Section~\ref{statement}.
First we discuss the notion of ``sequence of small points''
on semiabelian varieties,
and describe the major ingredients from which our result will follow.

\section{Small points on semiabelian varieties}
\label{smallpoints}

For the rest of this paper, $k$ denotes a number field.
Let $A$ be a semiabelian variety over $k$.
When $A$ is almost split, then after enlarging $k$ we have
an isogeny $\phi=(\phi_1,\phi_2): A \rightarrow A_0 \times \G_m^n$,
and we can define a canonical height $h: A(\kbar) \rightarrow \R$
by $h(x)=h_1(\phi_1(x)) + h_2(\phi_2(x))$,
where $h_1: A_0(\kbar) \rightarrow \R$ is the N\'eron-Tate canonical height
associated to a symmetric ample line bundle on $A_0$,
and $h_2: \G_m^n(\kbar) \rightarrow \R$
is the sum of the naive heights of the coordinates.
Then a sequence of points $x_i \in A(\kbar)$
is said to be a sequence of small points if $h(x_i) \rightarrow 0$.

But for more general semiabelian varieties over $k$,
there seems to be no natural canonical height with reasonable properties.
Therefore we instead use the definition of sequence of small points
given in~\cite{mlb}, which was based on a suggestion of Hrushovski.
For the details of the following, see~\cite{mlb}.

Embed $A$ as a quasi-projective variety in some $\PP^n$,
and let $h$ be a naive Weil height.
Fix an integer $m \ge 2$, and
let $[m]:A \rightarrow A$ denote multiplication by $m$.
Then there exists an integer $r \ge 1$ and real constants $M>0$ and $c>1$
such that for $z \in A(\kbar)$, $h(z)>M$ implies $h([m^r]z) > c h(z)$.
If $z \in A(\kbar)$, let $N(z)$ be the smallest integer $N \ge 1$
such that $h([m^N]z)>M$, or $\infty$ if no such $N$ exists.
We say that $\{z_i\}_{i \ge 1}$ is a {\em sequence of small points}
if $N(z_i) \rightarrow \infty$ in $\PP^1(\R)$.

It is shown in~\cite{mlb} that this notion is independent of choices.
Moreover, when $A$ is almost split, this notion agrees with the
one defined using canonical heights.

\section{The Mordell-Lang conjecture}

By the ``Mordell-Lang conjecture'' we mean the theorem below,
which was proved by McQuillan~\cite{mcquillan}, following work by
Faltings, Vojta, Hindry, Raynaud, and many others.
The division group $\Gamma'$ of a subgroup $\Gamma$ of a group $G$
is defined as
	$$\Gamma':= \{\, x \in G \mid
	\text{ there exists $n \ge 1$ such that } nx \in \Gamma \,\}.$$

\begin{theorem}[Mordell-Lang conjecture]
\label{mordelllang}
Let $A$ be a semiabelian variety over a number field $k$.
Let $\Gamma$ be a finitely generated subgroup of $A(\kbar)$,
and let $\Gamma'$ be the division group.
Let $X$ be a geometrically integral closed subvariety of $A$
that is not equal to
the translate of a sub-semiabelian variety (over $\kbar$).
Then $X(\kbar) \cap \Gamma'$ is not Zariski dense in $X$.
\end{theorem}

\begin{rems}
$\nichts$

\begin{enumerate}
\item	We obtain an equivalent statement if
``translate of a sub-semiabelian variety''
is replaced by
``translate of a sub-semiabelian variety by a point in $\Gamma'$.''
\item	We obtain an equivalent statement if
``not Zariski dense in $X$''
is replaced by
``contained in a finite union $\bigcup Z_j$
where each $Z_j$ is a translate of sub-semiabelian variety of $A_\kbar$,
and $Z_j \subseteq X_\kbar$.''
(Here $A_\kbar$ denotes $A \times_k \kbar$, and so on.)
\item	The theorem is true for {\em any} field $k$ of characteristic~0:
specialization arguments let one reduce to the number field case.
\item	There are function field analogues: see~\cite{hrushovski}.
\end{enumerate}
\end{rems}

\noindent
We will need only the weaker statement obtained by replacing
$\Gamma'$ by $\Gamma$ in Theorem~\ref{mordelllang}.
This is sometimes called the ``Mordellic part,''
because the special case where $X \subset A$
is a curve of genus $\ge 2$ in its Jacobian
is equivalent to Mordell's conjecture about the finiteness of $X(k)$.

\section{The Bogomolov and equidistribution conjectures}

If $X$ is a geometrically integral variety,
then a sequence of points $z_i$ is said to be {\em generic}
in $X$ if $z_i \in X(\kbar)$ for all $i$, and the $z_i$ converge
to the generic point of $X_\kbar$.
The latter condition means that
each subvariety $Y$ of $X_\kbar$ other than $X_\kbar$ itself
contains only finitely many points in the sequence.

\begin{conj}[Bogomolov conjecture for semiabelian varieties]
\label{bogomolov}
Let $A$ be a semiabelian variety over a number field $k$.
Let $X$ be a geometrically integral closed subvariety of $A$.
Let $\{z_i\}$ be a sequence of small points of $A(\kbar)$,
generic in $X$.
Then $X_\kbar$ is the translate of a sub-semiabelian variety of $A_\kbar$
by a torsion point.
\end{conj}

Bogomolov's original conjecture was for a curve of genus $\ge 2$
in its Jacobian.
This case was proved by Ullmo~\cite{ullmo}.
For $A$ an abelian variety or a torus,
Conjecture~\ref{bogomolov}
was proved by Zhang in~\cite{zhangbogomolov} and~\cite{zhangtori},
respectively.
Recently, a proof for the case where $A$ is almost split
was announced by Chambert-Loir~\cite{chambertloir2}.

Recall that if $\mu_i$ for $i \ge 1$ and $\mu$ are probability measures
on a metric space $X$,
then one says that the $\mu_i$ {\em converge weakly} to $\mu$
if for every bounded continuous function $f$ on $X$,
$\lim_{i \rightarrow \infty} \int f\mu_i = \int f \mu$.
For the following, we fix an embedding $\sigma: \kbar \hookrightarrow \C$,
and let $A_\sigma$ denote the semiabelian variety over $\C$
obtained by base extension by $\sigma$.

\begin{conj}[Equidistribution conjecture]
\label{equidistribution}
Let $A$ be a semiabelian variety over a number field $k$.
Let $\{z_i\}$ be a sequence of small points of $A(\kbar)$,
generic in $A$.
Let $\mu_i$ be the uniform probability measure
on the finite set $\sigma(G_k(z_i)) \subset A_\sigma(\C)$.
Then the $\mu_i$ converge weakly to the normalized Haar
measure $\mu$ on the maximal compact subgroup $A^0$ of $A(\C)$.
\end{conj}

\begin{rem}
The subgroup $A^0$ also equals the closure of $A_\sigma(\C)_\tors$
in the complex topology.
It is all of $A(\C)$ is $A$ is abelian,
and it is a ``polydisc'' if $A$ is a torus.
\end{rem}

The cases where $A$ is an abelian variety or a torus
are proved in~\cite{suz} and~\cite{bilu}, respectively.
A proof for the case where $A$ is an almost split semiabelian variety
has been announced by Chambert-Loir~\cite{chambertloir2}.
Bilu suggested in~\cite{bilu}
that Haar measure on $A^0$ should be the limit measure,
but he could not formulate Conjecture~\ref{equidistribution} precisely,
because at the time there was no notion
of ``sequence of small points'' for general semiabelian varieties.

\section{Statement of the result}
\label{statement}

The following was conjectured in the final section of~\cite{mlb}:

\begin{conj}[Mordell-Lang + Bogomolov]
\label{mainconj}
Let $A$ be a semiabelian variety over a number field $k$.
Let $\Gamma$ be a finitely generated subgroup of $A(\kbar)$,
and let $\Gamma'$ be its division group.
Let $X$ be a geometrically integral closed subvariety of $A$.
For $i \ge 1$, suppose that $x_i=\gamma_i+z_i \in X(\kbar)$
where $\gamma_i \in \Gamma'$
and $\{z_i\}_{i \ge 1}$ is a sequence of small points in $A(\kbar)$.
If $X_\kbar$ is not a translate of a sub-semiabelian variety of $A_\kbar$
by an element of $\Gamma'$,
then the $x_i$ are not Zariski dense in $X$.
\end{conj}

\begin{rem}
$\nichts$
\begin{enumerate}
\item	Equivalently, we could let
	$B_\epsilon:=\{\,z \in A(\kbar) \mid N(z)>1/\epsilon \,\}$
	and conjecture that for some $\epsilon>0$,
	$X(\kbar) \cap (\Gamma' + B_\epsilon)$
	is not Zariski dense in $X$.
	(Here $N$ is as in Section~\ref{smallpoints}.)
\item	It is equivalent if we remove the restriction that
	$X_\kbar$ not be the translate of a sub-semiabelian variety
	by an element of $\Gamma'$ and replace the conclusion
	``not Zariski dense in $X$''
	by ``contained in a finite union $\bigcup Z_j$
	where each $Z_j$ is a translate
	of a sub-semiabelian variety of $A_\kbar$
	by an element of $\Gamma'$,
	and $Z_j \subseteq X_\kbar$.''
\end{enumerate}
\end{rem}

Let $\fancyS(A)$ be the set of semiabelian varieties over number fields
that can be obtained from $A$
by taking algebraic subgroups, quotients, and products,
and changing the field of definition.

\begin{theorem}
\label{main}
Let $A$ be a semiabelian variety over a number field $k$.
Assume that Conjectures~\ref{bogomolov} and~\ref{equidistribution}
hold for all $B \in \fancyS(A)$.
Then Conjecture~\ref{mainconj} holds for subvarieties $X$ in $A$.
\end{theorem}

\begin{cor}
\label{maincor}
Conjecture~\ref{mainconj} holds when $A$ is almost split.
\end{cor}

\begin{proof}
One checks that all $B \in \fancyS(A)$ are almost split.
Chambert-Loir~\cite{chambertloir2} proved
Conjectures~\ref{bogomolov} and~\ref{equidistribution}
for almost split semiabelian varieties.
\end{proof}

\section{Properties of small points}

In preparation for the proof of Theorem~\ref{main},
we derive a few more properties of sequences of small points.
Throughout this section, $A$ denotes a semiabelian variety
over a number field $k$.
Let $\beta_n: A^n \rightarrow A^{n-1}$ be the map sending
$(x_1,\ldots,x_n)$ to $(x_2-x_1,x_3-x_1,\ldots,x_n-x_1)$.
Up to an automorphism of $A^{n-1}$,
this is the same as the map $\alpha_n$ used in~\cite{zhangbogomolov}.

Let $G_k=\Gal(\kbar/k)$.
For $x \in A(\kbar)$,
let $G_k(x)$ denote the $G_k$-orbit of $x$,
and let $\Diff_n(x)$ be the (arbitrarily ordered) list
of $[k(x):k]^n$ elements of $A^{n-1}$
obtained by applying $\beta_n$ to the elements of $G_k(x)^n$.
Given any sequence $x_1,x_2,\ldots$ of points in $A(\kbar)$,
let $\D_n=\D_n(\{x_i\})$
denote the infinite sequence obtained by concatenating
$\Diff_n(x_1)$, $\Diff_n(x_2)$, \dots.

\begin{lemma}
\label{sequences}
Let $A$ and $B$ be semiabelian varieties over a number field $k$.
Let $\{x_i\}$ be a sequence of small points in $A$,
and let $\{y_i\}$ be a sequence of small points in $B$.
\begin{itemize}
\item[(1)]	If $\{z_j\}$ is an infinite subsequence of $\{x_i\}$,
	then $\{z_j\}$ is a sequence of small points in $A$.
\item[(2)]	If $\{z_j\}$ is an infinite sequence
	obtained from $\{x_i\}$ by replacing each $x_i$
	with a finite number of copies of $x_i$,
	then $\{z_j\}$ is a sequence of small points in $A$.
\item[(3)]	If $\sigma \in G_k$, then $\{\sigma x_i\}$
	is a sequence of small points in $A$.
\item[(4)]	The sequence $\{(x_i,y_i)\}$
	is a sequence of small points in $A \times B$.
\item[(5)]	If $f: A \rightarrow B$ is a homomorphism,
	then $\{f(x_i)\}$ is a sequence of small points in $B$.
\item[(6)]	If $A=B$,
	then $\{x_i+y_i\}$ and $\{x_i-y_i\}$ are sequences of small points
	in $A$.
\item[(7)]	For any $n \ge 2$, $\D_n=\D_n(\{x_i\})$
	is a sequence of small points in $A$.
\item[(8)]	If $[k(x_i):k]$ is bounded, then there is a finite subset
	$T \subset A(\kbar)_\tors$ containing all but finitely many
	of the $x_i$.
\item[(9)]	There exists a sequence of positive integers $\{n_i\}$
		with $n_i \rightarrow \infty$,
		such that if $1 \le m_i \le n_i$,
		then the sequence $\{[m_i]x_i\}$ is still a sequence
		of small points.
\end{itemize}
\end{lemma}

\begin{proof}
Properties~(1) through~(4) are immediate from the definition,
and~(5) is proved in~\cite{mlb}.
Property~(6) follows from~(4) and~(5),
and~(7) follows from~(2), (3), (4), and~(6).

Next we prove~(8).
Let $S = \{\, x \in A(\kbar): h(x) \le M \,\}$.
By Northcott's Theorem, $s:=\#S$ is finite.
By definition of $N(x_i)$, $[m^j]x_i \in S$ for $j<N(x_i)$.
We have $N(x_i) > 1 + s$ for all but finitely many $i$.
For these $i$,
the pigeonhole principle yields $[m^j]x_i=[m^{j'}]x_i$ for some $j<j'\le 1+s$.
In particular $x_i \in T$,
where $T := \bigcup_{q < m^{1+s}} A(\kbar)[q]$.

Finally, we prove~(9).
Define $a_{ij}:=N([j]x_i)$.
By part~(5) applied to $[j]:A \rightarrow A$,
we have $\lim_{i \rightarrow \infty} a_{ij} = \infty$ for each $j$.
It follows formally that there is a sequence of positive integers $\{n_i\}$
tending to infinity,
such that $a_{i,m_i} \rightarrow \infty$ whenever $1 \le m_i \le n_i$.
\end{proof}

Next we have a sequence of lemmas leading up to Lemma~\ref{rationalpoints2},
which is the only other result from this section that will be used later.
Recall from~\cite{mlb} that one possible choice of $h$ is as follows.
We have an exact sequence
$0 \rightarrow T \rightarrow A \rightarrow A_0 \rightarrow 0$
where $T$ is a torus and $A_0$ is an abelian variety.
Enlarging $k$, we may assume that $T \isom \G_m^r$.
There is a compactification $\Abar$ to which $[2]$ extends,
equipped with effective line bundles $\LL_0$, $\LL_1$
with $\LL:=\LL_0 \tensor \LL_1$ ample,
such that $[2]^\ast \LL_0 = \LL_0^{\tensor 2}$
and $[2]^\ast \LL_1 = \LL_1^{\tensor 4}$.
In fact, $\Abar$ is also equipped with a map $\pibar: \Abar \rightarrow A_0$
extending $A \rightarrow A_0$,
and $\LL_1=\pibar^\ast \MM$,
where $\MM$ is a symmetric ample line bundle on $A_0$.
(See Section~1.1 of~\cite{mcquillan}, for example.)
Define $h=h_0+h_1$ where $h_0$ is a Weil height associated to $\LL_0$,
and $h_1$ is the {\em canonical height} associated to $\LL_1$
(i.e., the pullback of the canonical height on $A_0$ associated to $\MM$).

\begin{lemma}
\label{2x+a}
Given $a \in A(\kbar)$, there exists $M_a>0$ such that
if $x \in A(\kbar)$ and $h(x)>M_a$, then $h([2]x+a) > (3/2) h(x)$.
\end{lemma}

\begin{proof}
Translation-by-$a$ extends to a morphism $\tau_a:\Abar \rightarrow \Abar$,
and $\tau_a^\ast \LL_0 = \LL_0$.
It follows that $h_0(x+a)=h_0(x)+O(1)$,
where the $O(1)$ depends on $a$.
Also $[2]^\ast \LL_0 = \LL_0^{\tensor 2}$, so $h_0([2]x)=2 h_0(x) + O(1)$.
On the other hand, since $h_1$ is a quadratic function,
$h_1(x+a)=h_1(x)+O(h_1(x)^{1/2}) + O(1)$.
Hence
\begin{align*}
	h([2]x+a)	&= h_0([2]x+a) + h_1([2]x+a)	\\
			&= \left[ 2 h_0(x)+O(1) \right]
			+ \left[ 4 h_1(x) + O(h_1(x)^{1/2}) + O(1) \right] \\
			&= 2 h(x) +
			\left[ 2 h_1(x) + O(h_1(x)^{1/2}) + O(1) \right] \\
			&\ge 2 h(x) + O(1).
\end{align*}
\end{proof}

\begin{lemma}
\label{realvectorspace}
Let $\Gamma$ be a finitely generated subgroup of $A(\kbar)$,
and let $\{x_i\}$ be a sequence in $\Gamma'$.
If the image of $\{x_i\}$ in $\Gamma' \tensor \R = \Gamma \tensor \R$
converges to zero in the usual real vector space topology,
then $\{x_i\}$ is a sequence of small points.
\end{lemma}

\begin{proof}
Let $S:=\{\gamma_1,\gamma_2,\ldots,\gamma_n\} \subset \Gamma$
be a $\Z$-basis for $\Gamma/\Gamma_\tors$.
Let $U = \{\, \sum \epsilon_i \gamma_i : \epsilon_i \in \{-1,0,1\} \,\}$.
Let $f_1$, $f_2$, \dots, $f_u$ be the maps $A \rightarrow A$
of the form $x \mapsto [2]x+a$ for $a \in U$.
By repeated application of Lemma~\ref{2x+a},
we can find $M>0$ such that $h(x)>M$ implies $h(f_i(x))>(3/2)h(x)$
for all $i$.

Let $B$ be the subset of elements of $\Gamma'$ whose image
in $\Gamma \tensor \R$ have coordinates (with respect to the basis $S$)
bounded by~1 in absolute value.
For every $b_0 \in B$, there exists $i$ such that $b_1:=f_i(b_0) \in B$,
and then there exists $j$ such that $b_2:=f_j(b_1) \in B$,
and so on.
The intersection $I$ of $B$ with the finitely generated subgroup
generated by $b_0$ and $S$ is finite,
and $b_i \in I$ for all $i$.
But if $h(b_0)>M$, then $h(b_{m+1})>(3/2)h(b_m)$ for all $m$,
and in particular, the $b_i$ would be all distinct.
Hence $h(b_0) \le M$ for all $b_0 \in B$.
The lemma now follows from the definition of sequence of small points.
\end{proof}

\begin{rem}
The converse to Lemma~\ref{realvectorspace} is true, but we do not need it.
\end{rem}

\begin{lemma}
\label{rationalpoints1}
Let $\Gamma$ be a finitely generated subgroup of $A(\kbar)$.
Then $A(k) \cap \Gamma'$ is a finitely generated group.
\end{lemma}

\begin{proof}
By the Mordell-Weil theorem, $A_0(k)$ is finitely generated.
Since we have assumed that the toric part $T$ of $A$ is split,
Hilbert's Theorem~90 guarantees that there is no obstruction
to lifting generators of $A_0(k)$ to elements of $A(k)$.
Without loss of generality, enlarge $\Gamma$ to contain these lifts.
Then $A(k) \subseteq T(k)+\Gamma$,
and $A(k) \cap \Gamma' \subseteq (T(k) \cap \Gamma') + \Gamma$,
so we reduce to the case where $A = T = \G_m^r$.

Let $\pi_i:\G_m^r \rightarrow \G_m$ be the $i$-th projection.
Enlarging $\Gamma$ by replacing it by $\prod_{i=1}^r \pi_i(\Gamma)$,
we reduce to the case $r=1$; i.e., $A=\G_m$.
The finite generation of the unit group
and the finiteness of the class group of $k$
imply that $k^\ast$ is isomorphic as abstract group to
the direct sum of a finite torsion group
with a free abelian group of countable rank.
Any finite rank subgroup of such a group is finitely generated.
\end{proof}

\begin{lemma}
\label{rationalpoints2}
Let $\Gamma$ be a finitely generated subgroup of $A(\kbar)$.
Suppose $\{x_i\}$ is a sequence in $A(k)$,
and $x_i=\gamma_i + z_i$ where $\gamma_i \in \Gamma'$,
and $\{z_i\}$ is a sequence of small points in $A(\kbar)$.
Then there is a finitely generated subgroup of $\Gamma'$
containing all but finitely many of the $x_i$.
\end{lemma}

\begin{proof}
We may enlarge $k$ to assume $\Gamma \subset A(k)$.
Choose $\{n_i\}$ for $\{z_i\}$ as in part~(9) of Lemma~\ref{sequences}.
By elementary diophantine approximation (the pigeonhole principle),
there exist integers $m_i$ with $1 \le m_i \le n_i$,
and $\nu_i \in \Gamma$
such that the images of $m_i \gamma_i - \nu_i$ in $\Gamma \tensor \R$
approach zero as $i \rightarrow \infty$.
Then $\{ m_i \gamma_i -\nu_i \}$ is a sequence of small points
by Lemma~\ref{realvectorspace},
but $\{m_i z_i\}$ also is a sequence of small points,
so by part~(6) of Lemma~\ref{sequences},
$\{m_i x_i - \nu_i\}$ is a sequence of small points.
On the other hand, $m_i x_i -\nu_i \in A(k)$,
so by part~(8) of Lemma~\ref{sequences},
$x_i \in \Gamma'$ for all but finitely many $i$.
Finally, Lemma~\ref{rationalpoints1} implies that
there is a finitely generated subgroup of $\Gamma'$
containing all but finitely many $x_i$.
\end{proof}

\section{Measure-theoretic lemmas}
\label{measurelemmas}

We recall the following lemma from~\cite{mlb}:

\begin{lemma}
\label{measures}
Let $V$ be a projective variety over $\C$.
Let $S$ be a connected quasi-projective variety over $\C$.
Let $\Y \rightarrow V \times S$ be a closed immersion of $S$-varieties,
where $\Y \rightarrow S$ is flat with $d$-dimensional fibers.
For $i \ge 1$, let $s_i \in S(\C)$ and let $Y_i \subset V$ be
the fiber of $\Y \rightarrow S$ above $s_i$.
Let $\mu_i$ be a probability measure supported on $Y_i(\C)$.
If the $\mu_i$ converge weakly to a probability measure $\mu$ on $V(\C)$,
then the support of $\mu$ is contained in a $d$-dimensional Zariski closed
subvariety of $V$.
\end{lemma}

\begin{lemma}
\label{measures2}
Let $V$ and $S$ be quasi-projective varieties over $\C$, with $S$ integral.
Let $\Y$ be a subvariety of $V \times S$.
Let $s_1,s_2,\ldots$ be a sequence in $S(\C)$, Zariski dense in $S$.
Let $\mu_i$ be a probability measure with support contained in the fiber
of $\pi: \Y \rightarrow S$ above $s_i$, considered as subvariety of $V$.
Suppose the $\mu_i$ converge weakly to a probability measure $\mu$ on $V(\C)$.
Then the support of $\mu$ is contained in a subvariety of $V$
of dimension $\dim \Y - \dim S$.
\end{lemma}

\begin{proof}
Choose an embedding $V \hookrightarrow \PP^m$.
Without loss of generality, we may replace $V$ and $\Y$ by
their closures in $\PP^m$ and $\PP^m \times S$, respectively.
Replacing $S$ by a dense open subset $U$
and $\Y$ by $\pi^{-1}(U)$,
and passing to a subsequence,
we may assume that $\Y \rightarrow S$ is flat.
The result now follows from Lemma~\ref{measures}.
\end{proof}

\begin{lemma}
\label{measures3}
Retain the assumptions of the previous lemma, but assume in addition
that $V$ is a semiabelian variety over $\C$,
and that $\Y$ is not Zariski dense in $V \times S$.
Then $\mu$ does not equal the normalized Haar measure on
the maximal compact subgroup $V^0$ of $V(\C)$.
\end{lemma}

\begin{proof}
We have $\dim \Y - \dim S < \dim V$, so by the previous lemma,
$\mu$ is supported on a subvariety of $V$ of positive codimension.
But $V^0$ is Zariski dense in $V$.
\end{proof}

\section{Proof of Theorem~\ref{main}}

The proof will proceed through various reductions;
to aid the reader, we \fbox{box} cumulative assumptions
and other partial results to be used later in proof.
Let $G$ be the group of translations preserving $X$;
i.e., the largest algebraic subgroup of $A$ such that $X+G=X$.
We may assume \fbox{$\dim G=0$},
since otherwise we consider $X/G \hookrightarrow A/G$
and use part~(5) of Lemma~\ref{sequences}.
We may also enlarge $k$ to assume that \fbox{$\Gamma \subset A(k)$}.

If the theorem is false,
then there exists a sequence $x_i = \gamma_i + z_i \in X(\kbar)$,
generic in $X$,
with $\gamma_i \in \Gamma'$
and with $\{z_i\}$ a sequence of small points in $A(\kbar)$.
For $\sigma, \tau \in G_k$,
	$$\sigma x_i - \tau x_i = (\sigma \gamma_i - \tau \gamma_i)
				+ (\sigma z_i - \tau z_i).$$
Some multiple of $\gamma_i$ is in $\Gamma \subset A(k)$,
so $\sigma \gamma_i - \tau \gamma_i$ is torsion.
Applying part~(7) of Lemma~\ref{sequences} to $\D_2(\{z_i\})$,
and then applying parts~(2) and~(6),
we find that $\D_2:=\D_2(\{x_i\})$ is a sequence of small points in $A$.
Then by parts~(2) and~(4) of Lemma~\ref{sequences},
$\D_n$ is a sequence of small points for each $n \ge 2$.
Fix \fbox{$n>\dim A$}.

By repeated application of Conjecture~\ref{bogomolov},
we may discard finitely many of the $x_i$ in order to assume
that the Zariski closure $\Dbar_n$ of $\D_n$ in $A_\kbar$ is a finite union
$\bigcup_{j=1}^s (B_j +t_j)$
where $B_j$ is a sub-semiabelian variety of $A_\kbar$,
and $t_j \in A(\kbar)$ is a torsion point.
If we replace $X$ by the image of $X$ under multiplication by
a positive integer $N$, and replace each $x_i$ by $Nx_i$,
then $\bigcup_{j=1}^s B_j$ is unchanged.
If we pass to a subsequence of the $x_i$ or enlarge $k$,
then the new $\bigcup_{j=1}^s B_j$ can only be smaller.
Since $A_\kbar$ is noetherian, we may assume
without loss of generality that these
operations are done so as to make $\bigcup_{j=1}^s B_j$ minimal.
Moreover, by multiplying by a further integer $N$
we may assume that $t_j=0$ for all $j$.
Now, any further operations of the types above will leave $\Dbar_n$ unchanged,
equal to $\bigcup_{j=1}^s B_j$.
Enlarging $k$, we may assume that each $B_j$ is defined over $k$.

Repeating the same procedure with $2$ instead of $n$,
we may minimize $\Dbar_2 = \bigcup_{j=1}^u C_j$,
and assume that each $C_j$ is a sub-semiabelian variety of $A$.
We now show that there is only one $C_j$.
By the pigeonhole principle, there is some $C_j$, say $C:=C_1$,
such that for infinitely many $i$,
at least a fraction $1/u$ of the elements of $G_k(x_i)^2$ are mapped
by $\beta_2$ into $C$.
Passing to a subsequence of the $x_i$,
we may suppose that this holds for {\em all} $i$.
Let $\pi:A \rightarrow A/C$ be the projection, and let $y_i=\pi(x_i)$.
If $\ell$ is a finite Galois extension of $k$ containing $k(x_i)$,
then it follows that $\sigma y_i - \tau y_i = 0$
for at least a fraction $1/u$ of the pairs $(\sigma,\tau) \in \Gal(\ell/k)^2$.
Thus the subgroup of $\Gal(\ell/k)$ stabilizing $y_i$ must have index
at most $u$.
Hence $\deg_k(y_i) = \# G_k(y_i) \le u$.
Then $\D_2(\{y_i\})$ consists of points of degree bounded by $u^2$.
On the other hand, $\D_2(\{y_i\})$ is a sequence of small points in $A/C$,
by parts~(5) and~(7) of Lemma~\ref{sequences}.
By part~(8) of Lemma~\ref{sequences},
there is a finite subset $T$ of torsion points of $A/C$
such that $\Diff_2(y_i) \subseteq T$ for all sufficiently large $i$.
Passing to a subsequence of the $x_i$,
and multiplying everything by an integer $N$ to kill $T$,
we may assume that $\Diff_2(y_i)=\{0\}$ for all $i$.
Then $\Diff_2(x_i) \subseteq C$, so $\Dbar_2 \subseteq C$,
and hence \fbox{$\Dbar_2=C$} by definition of $C$.

We next show that there is only one $B_j$, and that it equals $C^{n-1}$.
By definition of $\D_2$ and $\D_n$, we have $B_j \subseteq C^{n-1}$
for each $j$.
By the pigeonhole principle, there is some $B_j$, say $B:=B_1$,
such that for infinitely many $i$,
at least a fraction $1/s$ of the elements of $G_k(x_i)^n$
are mapped by $\beta_n$ into $B$.
Passing to a subsequence, we may suppose that this holds for {\em all} $i$.
For $1 \le q \le n-1$, define the ``coordinate axis''
$C_{(q)} = 0 \times \cdots \times 0 \times C \times 0 \times \cdots \times 0$,
with $C$ in the $q$-th place.
Let $B_{(q)}=B \cap C_{(q)}$.
By the pigeonhole principle again,
given $i$,
there exist $w_1,w_2,\ldots,w_q,w_{q+2},\ldots,w_n \in G_k(x_i)$
such that $\beta_n(w_1,w_2,\ldots,w_q,\zeta,w_{q+2},\ldots,w_n) \in B$
for at least a fraction $1/s$ of the elements $\zeta$ of $G_k(x_i)$.
Subtracting,
we find that $\beta_n(0,0,\ldots,0,\zeta-\zeta',0,\ldots,0) \in B_{(q)}$
for at least a fraction $1/s^2$
of the pairs of elements $\zeta,\zeta'$ of $G_k(x_i)$.
As before, this implies
(after passing to a subsequence and multiplying by a positive integer again)
that the image $y_i$ of $x_i$ in $C/B_{(q)}$ satisfies $\Diff_2(y_i)=\{0\}$.
Then $\Dbar_2 \subseteq B_{(q)} \subseteq C$, so $B_{(q)}=C$.
This holds for all $q$, so \fbox{$\Dbar_n=B=C^{n-1}$}.

If $C=\{0\}$, then $\D_2=\{0\}$, and then by definition of $\D_2$,
$x_i \in X(k)$ for all $i$.
Lemma~\ref{rationalpoints2} implies that
all but finitely many $x_i$
are contained in a finitely generated subgroup $\tilde{\Gamma}$ of $\Gamma'$.
The Mordellic part of Theorem~\ref{mordelllang} applied to $\tilde{\Gamma}$
implies that $X_\kbar$ is a translate of a sub-semiabelian variety.
But $\dim G=0$, so $X$ is a point.
Moreover this point is in $\Gamma'$
(by our application of Lemma~\ref{rationalpoints2}),
so we contradict the hypotheses.
Therefore we may assume \fbox{$\dim C \ge 1$}.

Let $S=\pi(X) \subseteq A/C$.
Note that $S$ is integral.
Consider the fibered power $X^n_S := X \times_S X \times_S \cdots \times_S X$
as a subvariety of $X^n$.

Let \fbox{$m=\dim C$}.
Note that $1 \le m \le \dim A < n$.
Let $\dim(X/S)$ denote the relative dimension;
i.e., the dimension of the generic fiber of $X \rightarrow S$.
Then \fbox{$\dim(X/S)<m$}, since otherwise $X$ (being closed)
would equal the entire inverse image of $S$ under $A \rightarrow A/C$,
and then $C \subseteq G$, contradicting $\dim G=0$.
Hence
\begin{equation}
\label{diminequality}
	\dim(X^n_S/S) = n \dim(X/S) \le n(m-1) < m(n-1)
	= \dim(C^{n-1} \times S/S).
\end{equation}
The homomorphism $\beta_n$ restricts
to a morphism $X^n_S \rightarrow C^{n-1}$.
We also have the obvious morphism $X^n_S \rightarrow S$.
Let $\Y$ denote the image of the product morphism
$X^n_S \rightarrow C^{n-1} \times S$.
Then $\dim \Y < \dim C^{n-1} \times S$, by~(\ref{diminequality}).

Since $\sigma x_i - \tau x_i \in C(\kbar)$ for all $\sigma,\tau \in G_k$
and all $i \ge 1$,
concatenating the finite subsets
$\beta(G_k(x_i)^n) \times \{\pi(x_i)\}$ of $C^{n-1} \times S$
yields a sequence of points $y_j=(c_j,s_j)$ in $\Y$.
The $c$-sequence is simply $\D_n$.
and each $s_j$ equals $\pi(x_i)$ for some $i$.
Now fix an embedding $\sigma:\kbar \hookrightarrow \C$,
and let $\mu_j$ be the uniform probability measure on
the finite subset $\sigma(G_k(c_j)) \subset C_\sigma^{n-1}(\C)$.
Conjecture~\ref{equidistribution} implies that the $\mu_j$
converge to the normalized Haar measure $\mu$
on the maximal compact subgroup of $C_\sigma^{n-1}(\C)$.

On the other hand, $\mu_j$ is supported on $\sigma\beta(G_k(x_i)^n)$,
which is contained in the fiber of $\Y \rightarrow S$ above $s_j$,
when we consider the fiber as a subvariety of $V:=C_\sigma^{n-1}$.
Lemma~\ref{measures3} implies that the $\mu_j$ cannot converge
to $\mu$.
This contradiction completes the proof of Theorem~\ref{main}.


\end{document}